\documentclass[reqno,12pt]{amsart}

\title[Hurwitz-Hodge integrals, $E_{6}$, $D_{4}$, and the
CRC]{Hurwitz-Hodge integrals, the $E_{6}$ and $D_{4}$ root systems,
and the Crepant Resolution Conjecture}

\author{Jim Bryan and Amin
Gholampour}
\address{
Dept of Math, Univ. of British Columbia,
Vancouver, BC, Canada 
}
\email{jbryan@math.ubc.ca}
\email{amin@math.ubc.ca}

\usepackage{amsmath,amsthm,amsfonts,amscd}
\usepackage{eepic}
\usepackage{xy}
\usepackage{graphicx}
\usepackage[small]{diagrams}
\usepackage{verbatim}
\usepackage{times}

\newtheorem{thm}{Theorem}[section]
\newtheorem{theorem}{Theorem}
\newtheorem{lem}[thm]{Lemma}
\newtheorem{prop}[thm]{Proposition}

\newtheorem{lemma}[theorem]{Lemma}
\newtheorem{corollary}[theorem]{Corollary}
\theoremstyle{definition}
\newtheorem{rem}[thm]{Remark}

\newtheorem{remark}[theorem]{Remark}

\newcommand{\X}{\mathcal{X}}
\newcommand{\s}{\sigma}
\newcommand{\T}{\tau}
\newcommand{\z}{\zeta}
\newcommand{\R}{\rho}
\newcommand{\C}{\mathbb{C}}
\newcommand{\Z}{\mathbb{Z}}

\newcommand{\ev}{\operatorname{ev}}
\renewcommand{\H}{\overline{H}}
\newcommand{\M}{\overline{M}}
\renewcommand{\P}{\mathbb{P}^{1}}
\renewcommand{\hat}{\widehat}
\newcommand{\Cbar}{\overline{C}}
\newcommand{\ZtwoZtwo}{\Z_{2} \times \Z _{2}}

\newcommand{\GHilb}{\operatorname{G-Hilb}}
\newcommand{\orb}{\operatorname{orb}}
\newcommand{\LL}{\mathsf{Length}}
\newcommand{\LangRang}[1]{\left\langle #1 \right\rangle}
\newcommand{\Evee }{\mathbb{E}^{\vee }}

\newcommand{\jbar}{\overline{\jmath}}
\newcommand{\h}{\mathbf{h}}
\newcommand{\leftbig}[2]{\left#1\rule{0in}{#2pt}\right.}
\newcommand{\rightbig}[2]{\left.\rule{0in}{#2pt}\right#1}

\newcommand{\pimapA}{%
\begin{picture}(0,0)(-8,-6)
\qbezier(10,-37)(0, -15)(0, 0)
\qbezier(10,37)(0, 15)(0, 0)
\put(-15,3){\makebox(0,0)[t]{$\overline{\pi }$}}
\put(12.3,-38.5){\vector(1,-1){0}}
\end{picture}
}

\begin{document}
\maketitle 
\begin{abstract}
Let $G$ be the group $A_{4}$ or $\ZtwoZtwo $. We compute the integral
of $\lambda _{g}$ on the Hurwitz locus $\H _{G}\subset \M_{g}$ of
curves admitting a degree 4 cover of $\P $ having monodromy group
$G$. We compute the generating functions for these integrals and write
them as a trigonometric expression summed over the positive roots of
the $E_{6}$ and $D_{4}$ root systems respectively. As an application,
we prove the Crepant Resolution Conjecture for the orbifolds $[\C
^{3}/A_{4}]$ and $[\C ^{3}/ (\ZtwoZtwo )]$.
\end{abstract}


\section{Introduction}

In his seminal 1983 paper \cite{Mumford}, Mumford developed an
enumerative geometry for the moduli space of curves analogous to
Schubert calculus. On the Grassmannian, one can integrate the Chern
classes of the tautological bundle over Schubert cycles, namely cycles
given by the loci of linear spaces satisfying various incidence
conditions. On the moduli space of stable curves, one can integrate
Chern classes of the Hodge bundle over Hurwitz cycles, namely cycles
defined by the loci of curves satisfying some Hurwitz conditions. Such
integrals (and their variants) are called Hurwitz-Hodge integrals and
they arise in various contexts, notably in orbifold Gromov-Witten
theory, e.g.  \cite{Bryan-Graber-Pandharipande,
Cavalieri-Hurwitz-Hodge,CCIT-twisted, CCIT-CRC}.

In this paper, we consider Hurwitz-Hodge integrals over a natural
class of $g$ dimensional cycles on $\M _{g}$ which are defined as
follows.

Let $G$ be either $A_{4}$, the alternating group on 4 letters, or
$\ZtwoZtwo $, the Klein four group. Let
\[
H_{G}\subset M_{g}
\]
be the locus of genus $g$ curves $C$
admitting a degree 4 map 
\[
f:C\to \P 
\]
whose monodromy group
is contained in $G$.

The branch points $p_{1},\dotsc ,p_{n}\in \P $ of $f$ are then such
that $f^{-1} (p_{i})$ consists of exactly 2 points. By the
Riemann-Hurwitz formula,
\[
g=n-3
\]
and consequently $H_{G}$ has dimension $g$. 

In modern terms,
$H_{G}$ can be described as $M_{0,n} (BG)$, the moduli stack
of twisted maps to the the classifying stack $BG$. As such,
there is a natural compactification $\H _{G}\subset \M _{g}$ given
by twisted stable maps $\M _{0,n} (BG)$.

The $G$-Hurwitz space 
\[
\H _{G}=\M _{0,n} (BG)
\]
has components indexed by the monodromy around the $n$ points. These
are given by $n$-tuples of non-trivial conjugacy classes in $G$. Since
each component has dimension $g$, we can evaluate the Hodge class
$(-1)^{g}\lambda _{g}$ on each component to obtain a rational
number. There are three non-trivial conjugacy classes in $G$, so the
natural generating functions for these $G$-Hurwitz-Hodge integrals are
formal power series $F_{G} (x_{1},x_{2},x_{3})$ in three variables
(defined in detail in \S\ref{sec notation and results}). Our
main result is an explicit formula for $F_{G}$ written in terms of the
$E_{6}$ and $D_{4}$ root systems for $G$ equal to $A_{4} $ and
$\ZtwoZtwo $ respectively.

To write our expression for $F_{G}$, we will need to introduce
some concepts which will relate conjugacy classes of $G$ to the
$E_{6}$ and $D_{4}$ root systems. 

Both $A_{4}$ and $\ZtwoZtwo $ are naturally subgroups of $SO (3)$
(they are the symmetry groups of the tetrahedron and the prism over
the 2-gon). Let $\hat{G}$ be the binary version of $G$, that is the
preimage of $G$ in $SU (2)$ (namely the binary tetrahedral group
$\hat{A}_{4}$ and the quaternion 8 group $Q$).
\[
\begin{diagram}
\hat{G}&\rInto & SU (2)\\
\dTo & &\dTo \\
G&\rInto & SO (3)
\end{diagram}
\]

By the classical McKay correspondence
\cite{Gonzalez-Sprinberg-Verdier,McKay,Reid-asterisque}, finite
subgroups of $SU (2)$ admit an ADE classification where the
non-trivial irreducible representations of a group naturally
correspond to the nodes of the associated Dynkin diagram. In this
classification, the binary tetrahedral group $\hat{A}_{4}$ and the
quaternion 8 group $Q$ correspond to the $E_{6}$ and $D_{4}$ Dynkin
diagrams respectively. The non-trivial irreducible representations of
$\hat{A}_{4}$ and $Q$ that pullback from representations of $A_{4}$
and $\ZtwoZtwo $ correspond to the white nodes in the Dynkin diagrams
below:

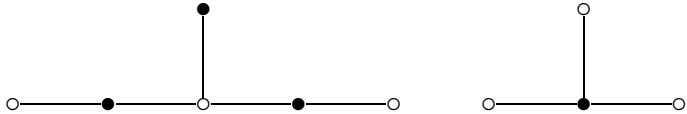
\vspace{-0.5cm}
\begin{center}
\begin{diagram}[height=1.5em,width=1.5em,abut]
      &       &        &       &\bullet&       &        &       &      &&      &        &\circ   & &     \\
      &       &        &       &\dLine &       &        &       &      &&      &        &\dLine  &&      \\
\circ &\hLine &\bullet &\hLine &\circ  &\hLine &\bullet &\hLine &\circ &&\circ &\hLine  &\bullet &\hLine &\circ \\
\end{diagram}
\end{center}

ADE Dynkin diagrams also correspond to simply laced root systems where
the nodes of the diagram correspond to simple roots of the root
system. Let $R$ be the $E_{6}$ or $D_{4}$ root system, let $\rho $ be
a non-trivial irreducible representation of $G$, and let $e_{\rho }$
be the simple root which corresponds to the same node in the Dynkin
diagram as $\rho $. For any positive root $\alpha \in R^{+}$, let
$\alpha ^{\rho }$ denote the coefficient of $e_{\rho }$ in $\alpha $.

We index non-trivial conjugacy classes of $G$ by $i\in \{1,2,3 \}$ and
we let $\chi _{\rho }^{i}$ be the value of the character of a
representation $\rho $ on the $i$th conjugacy class. Let $z_{i}$ be
the order of the centralizer of the $i$th conjugacy class, and let $V$
be the 3 dimensional representation of $G$ arising from the embedding
$G\subset SO (3)$. We define the following matrix which is a
modification of the character table of $G$:
\[
L_{\rho }^{i} = \frac{1}{z_{i}}\sqrt{3-\chi _{V}^{i}}\,\,\, \chi _{\rho }^{i}\,.
\]

Our main result is the following
\begin{theorem}\label{thm: main theorem root system formulation}
Let $G$ be $A_{4}$ or $\ZtwoZtwo $ and let $R$ be the $E_{6}$ or
$D_{4}$ root system respectively. The generating function for the
$G$-Hurwitz-Hodge integrals is given by
\[
F_{G} (x_{1},x_{2},x_{3}) = \frac{1}{2}\sum _{\alpha \in R^{+}}\h
\left(\pi +\sum _{\rho }\alpha ^{\rho }\left(\frac{2\pi \dim \rho
}{|G|}+ \sum _{i}L_{\rho }^{i}x_{i}\right)
\right)
\]
where $R^{+}$ is the set of positive roots of $R$, the sum over $\rho
$ is over non-trivial irreducible representations of $G$, and $\h (u)$
is the series defined by
\[
\h''' (u) = \frac{1}{2}\tan \left(\frac{-u}{2} \right).
\]
\end{theorem}
The above formula is expanded out explicitly for $\ZtwoZtwo $ in
Proposition~\ref{prop: explicit formula for FZ2Z2} and for $A_{4}$ in
Proposition~\ref{prop: explicit formula for FA4}. Note that since the
constant, linear, and quadratic terms of the series $\h (u)$ are
undefined, the same is true for $F_{G} (x_{1},x_{2},x_{3})$. This
corresponds to the fact that $\M _{0,n} (BG)$ is not defined for
$n<3$.

We prove in Proposition~\ref{prop: orbifold GW invs = Hurwitz-Hodge
ints} that $F_{G} (x_{1},x_{2},x_{3})$ is equal to the (non-classical
part of the) genus zero Gromov-Witten potential for the orbifold $[\C
^{3}/G]$. In \cite{Bryan-Graber} it is conjectured that for a crepant
resolution $Y\to X$ of an orbifold $\X $ satisfying the Hard Lefschetz
condition, the Gromov-Witten potentials of $Y$ and $\X $ are related
by a linear change of variables and a specialization of quantum
parameters of $Y$ to roots of unity.

The singular space $X=\C ^{3}/G$ underlying the orbifold $\X =[\C
^{3}/G] $ admits a preferred Calabi-Yau resolution
\[
\pi :\GHilb (\C ^{3}) \to X
\]
given by Nakamura's Hilbert scheme of $G$ clusters \cite{BKR}. In
\cite{Bryan-Gholampour3}, we completely compute the Gromov-Witten
theory of $\GHilb (\C ^{3})$ for all finite subgroups $G\subset SO
(3)$ and we use the Crepant Resolution Conjecture to obtain a
predicition for the genus zero Gromov-Witten potential of the orbifold
$\X =[\C ^{3}/G]$. The change of variables matrix for the conjecture
in this example is given by $\sqrt{-1}L_{\rho }^{i}$ and the roots of
unity are given by
\[
q_{\rho } = \exp \left(\frac{2\pi i\dim \rho }{|G|} \right).
\]
Note that the matrix $L_{\rho }^{i}$, the roots of unity $q_{\rho }$,
and the formula in Theorem~\ref{thm: main theorem root system
formulation} make sense for any finite subgroup $G\subset SO (3)$
(although the number of variables differs from three in
general). Indeed, the conjectural formula for $F_{\X } $, the
(non-classical part of the) Gromov-Witten potential of $\X =[\C
^{3}/G]$ is exactly the same as the formula for $F_{G}$ in
Theorem~\ref{thm: main theorem root system formulation} (see the
conjecture in \cite{Bryan-Gholampour3}).

An immediate consequence of Proposition~\ref{prop: orbifold GW invs =
Hurwitz-Hodge ints} and Theorem~\ref{thm: main theorem root system
formulation} is the proof of the conjecture in
\cite{Bryan-Gholampour3}. In particular, we have proven the following.
\begin{theorem}\label{thm: CRC is true for A4 and Z2xZ2}
The genus zero Crepant Resolution Conjecture is true for the orbifold
$[\C ^{3}/G] $ with its crepant resolution given by $\GHilb (\C ^{3})$
when $G$ is $\ZtwoZtwo $ or $A_{4}$.
\end{theorem}

The Gromov-Witten invariants for  $\X =[\C ^{3}/G]$
in general can also be described as integrals over $G$-Hurwitz loci in
$\M _{g}$ but for other $G\subset {SO (3)}$, the integral of $\lambda
_{g}$ is replaced with a slightly more exotic Hodge classes obtained
from Chern classes of eigen-subbundles of the Hodge bundle.

\section{Notation and Results}\label{sec notation and results} Let $G$
be a finite group, and $\overline{M}_{0,\,n}(BG)$ be the moduli space
of genus $0$, $n$-marked twisted stable maps to $BG$. The evaluation
maps, denoted by $\ev_{i}$ for $i\in \{1,\dotsc ,n \}$, take values in
the inertia stack $IBG$. The coarse space of this stack is a finite
collection of points, one for each conjugacy class in $G$.

Let $S=(c_1,\dots,c_n)$ be an $n$-tuple of  conjugacy classes in $G$.
We define the following open and closed substack of
$\overline{M}_{0,\,n}(BG)$:
\[
\overline{M}_{S}(BG)=\bigcap_{i=1}^n \; \ev_{i}^{-1}(c_i).
\]

Concretely, $\M _{S} (BG)$ parametrizes $G$ covers of an $n$ marked
genus zero curve with monodromy $c_{i}$ around the $i$th marked point.

In this paper we will deal with the cases that the group $G$ is either
$\ZtwoZtwo $, $A_4$, or $S_{4}$. We fix a notation for the conjugacy
classes in these groups that we will use throughout the paper:
\begin{itemize}
\item Let $1,\, \T, \, \s, \, \R, \, \z$, denote the conjugacy classes
in $S_4$ corresponding to the elements $(1), \; (1\, 2),\; (1\, 2\,
3), \; (1\, 2\, 3 \, 4), \; (1\, 2)(3\,4)$, respectively.
\item Let $1,\, \s_1, \, \s_2, \, \z$, denote the conjugacy classes in
$A_4$ corresponding to the elements $(1),\;(1\, 2\, 3),\;(1\, 3\,
2),\;(1\, 2)(3\,4)$, respectively.
\item Let $1,\, \z_1, \, \z_2, \, \z_3$, denote the conjugacy classes
in $\ZtwoZtwo $ corresponding to its four elements.
\end{itemize}

All the above groups have a natural action on the set of four
elements. Thus to any element
\[
[f:C \to BG]\in \M _{S} (BG)
\]
we can associate a degree 4 cover
\[
\Cbar \to C
\]
with monodromy type $c_i$ over the $i$-th marked point. Let
\[
\overline{\pi } :\overline{\mathcal{C}}\to \M _{S} (BG)
\]
be the universal family of the four fold covers and let
\[
\Evee  = R^{1}\overline{\pi } _{*} (\mathcal{O}_{\overline{\mathcal{C}}})
\]
be the dual Hodge bundle of this family.

We now define the \emph{$G$-Hurwitz-Hodge integrals} as follows:
\[
\LangRang{c_{1}\dotsb c_{n}}
^{G}=\int_{[\overline{M}_{S}(BG)]}
c (\Evee ).
\]
When $G$ is $A_{4}$ or $\ZtwoZtwo $ and $\overline{C}$ is connected,
it is genus $g$ and the integrand in the above definition is
$(-1)^{g}\lambda _{g}$.  The conjugacy classes $c_{i} $ in
$\LangRang{c_{1}\dotsb c_{n}}$ are called \emph{insertions} and the
total number of insertions will be called the \emph{length}. We will
drop the superscript $G$, when it is understood from context.

We define $F_{G}$, the generating functions for these integrals, by
\begin{align*}
F_{A_{4}} (x_{1},x_{2},x_{3})& = \sum _{n_{1}+n_{2}+n_{3}\geq 3}
\LangRang{ \sigma _{1}^{n_{1}}\sigma _{2}^{n_{2}}\zeta ^{n_{3}}
}^{A_{4}}\quad \frac{x_{1}^{n_{1}}}{n_{1}!}
\frac{x_{2}^{n_{2}}}{n_{2}!} \frac{x_{3}^{n_{3}}}{n_{3}!},\\
F_{\ZtwoZtwo } (x_{1},x_{2},x_{3})& = \sum _{n_{1}+n_{2}+n_{3}\geq 3}
\LangRang{ \zeta _{1}^{n_{1}}\zeta  _{2}^{n_{2}}\zeta_{3} ^{n_{3}}
}^{\ZtwoZtwo }\quad \frac{x_{1}^{n_{1}}}{n_{1}!}
\frac{x_{2}^{n_{2}}}{n_{2}!} \frac{x_{3}^{n_{3}}}{n_{3}!},\\
F_{S_{4}} (x_{1},x_{2},x_{3},x_{4})& = \sum _{n_{1}+n_{2}+n_{3}+n_{4}\geq 3}
\LangRang{ \tau ^{n_{1}}\sigma ^{n_{2}}\rho ^{n_{3}}\zeta ^{n_{4}}
}^{S_{4}}\quad \frac{x_{1}^{n_{1}}}{n_{1}!}
\frac{x_{2}^{n_{2}}}{n_{2}!} \frac{x_{3}^{n_{3}}}{n_{3}!} \frac{x_{4}^{n_{4}}}{n_{4}!}.
\end{align*}

Our use of $F_{S_{4}}$ is auxiliary to our computations of $F_{A_{4}} $
and $F_{\ZtwoZtwo }$ and we do not determine it completely.
 
For concreteness, we write out the formula in Theorem~\ref{thm: main
theorem root system formulation} explicitly for the two cases of
$A_{4}$ and $\ZtwoZtwo $.  As before, we define the series $\h (u)$ by
\[
\h''' (u) = \frac{1}{2}\tan \left(-\frac{u}{2} \right).
\]
By a theorem of Faber and Pandharipande
\cite{Faber-Pandharipande-logarithmic}, $\h (u)$ is the generating
series\footnote{Note that our series $\h (u)$ is equal to $-u^{2}H
(u)$ where $H (u)$ is the series defined in Faber-Pandharipande
\cite{Faber-Pandharipande-logarithmic}.} for $\Z_2$-Hurwitz-Hodge
integrals, namely the integral of $-\lambda _{g}\lambda _{g-1}$ over
the hyperelliptic locus $\H _{\Z_2}\subset \M _{g}$.

\begin{prop}\label{prop: explicit formula for FZ2Z2}
The generating function for $\ZtwoZtwo $-Hurwitz Hodge integrals is given by the formula
\begin{align*}
F_{\ZtwoZtwo }  =&\quad \h \left(\frac{1}{2} (x_{1}+x_{2}+x_{3})-\frac{\pi }{2} \right) + \h \left(\frac{1}{2} (-x_{1}+x_{2}-x_{3})-\frac{\pi }{2} \right) \\
&+ \h \left(\frac{1}{2} (x_{1}-x_{2}-x_{3})-\frac{\pi }{2} \right) + \h \left(\frac{1}{2} (-x_{1}-x_{2}+x_{3})-\frac{\pi }{2} \right) \\
&+\frac{1}{2} \h(x_{1})+\frac{1}{2} \h(x_{2})+\frac{1}{2} \h(x_{3}).
\end{align*}
\end{prop}

\newcommand{\xplusx}{\, \,\,\, x_1\,\,+\,\, x_2}
\begin{prop}\label{prop: explicit formula for FA4}
Let $\omega =e^{2\pi i/3}$. The generating function for
$A_{4}$-Hurwitz Hodge integrals is given by the formula
\begin{align*}
F_{A_{4}}=&\quad \h\left(\frac{1}{\sqrt{3}}(\xplusx )+\frac{x_3}{2}-\frac{5 \pi}{6}\right)\,\,
+2\h\left(\frac{1}{\sqrt{3}}(\!\xplusx )-\frac{\pi}{3}\right)\\
&+\h\left(\frac{1}{\sqrt{3}}(\omega x_1+\overline{\omega} x_2)+\frac{x_3}{2}-\frac{5 \pi}{6}\right)
+2\h\left(\frac{1}{\sqrt{3}}(\omega x_1+\overline{\omega} x_2)-\frac{\pi}{3}\right)\\
&+\h\left(\frac{1}{\sqrt{3}}(\overline{\omega} x_1+\omega x_2)+\frac{x_3}{2}-\frac{5 \pi}{6}\right)
+2\h\left(\frac{1}{\sqrt{3}}(\overline{\omega} x_1+\omega x_2)-\frac{\pi}{3}\right)\\
&+\h\left(\frac{1}{\sqrt{3}}(\xplusx )-\frac{x_3}{2}+\frac{\pi}{6}\right)\\
&+\h\left(\frac{1}{\sqrt{3}}(\omega x_1+\overline{\omega} x_2)-\frac{x_3}{2}+\frac{\pi}{6}\right)\\
&+\h\left(\frac{1}{\sqrt{3}}(\overline{\omega} x_1+\omega x_2)-\frac{x_3}{2}+\frac{\pi}{6}\right)
\,\,+4\h\left(\frac{x_3}{2}+\frac{\pi}{2}\right) +\frac{1}{2}\h\left(x_{3}\right).
\end{align*}
\end{prop}

\subsection{Outline of the proof} We prove Theorem~\ref{thm: main
theorem root system formulation} by first computing $F_{\ZtwoZtwo }$
to prove Proposition~\ref{prop: explicit formula for FZ2Z2} and then
by computing $F_{A_{4}}$ to prove Proposition~\ref{prop: explicit
formula for FA4}. For each of $F_{\ZtwoZtwo }$ and $F_{A_{4}}$ we
first prove that $F$ is uniquely determined by the WDVV equations
along with certain specializations. We then prove that our formulas
for $F$ satisfy the WDVV equations and specialize correctly. The
required specialization for $F_{\ZtwoZtwo }$ uses the
Faber-Pandharipande computation of $F_{\Z_2}$. The required
specializations for $F_{A_{4}}$ uses the previously proven formula for
$F_{\ZtwoZtwo }$ as well as certain generating series for
$S_{4}$-Hurwitz-Hodge integrals. These are in turn determined by a
WDVV argument in \S\ref{sec: computing S4} use the
$\Z_3$-Hurwitz-Hodge integrals computed in
\cite{Bryan-Graber-Pandharipande}.

\section{The WDVV equations}\label{sec: WDVV} The $G$-Hurwitz-Hodge
integrals $\LangRang{c_{1}\dotsb c_{n}}^{G}$ satisfy the following
version of the WDVV equations. These equations are the primary tools
of this paper.
\begin{thm}\label{thm: WDVV sum form}
For any $n$-tuple $(c_{1},\dotsc ,c_{n})$ of conjugacy classes of $G$
and any subset $I\subset \{1,\dotsc ,n \}$ of cardinality $|I|$, let
$c_{I}$ denote the corresponding $|I|$-tuple of conjugacy classes and
let $I^{c}$ be the complement of $I$. For $g\in G$ let $(g)$ denote
the corresponding conjugacy class and let $z (g)$ be the order of the
centralizer of $g$. Then
\[
\LangRang{c_{1}\dotsb c_{n} (a_{1}a_{2}|a_{3}a_{4})} ^{G} =\LangRang{c_{1}\dotsb c_{n} (a_{1}a_{3}|a_{2}a_{4})} ^{G} 
\]
where $\LangRang{c_{1}\dotsb c_{n} (a_{i}a_{j}|a_{k}a_{l})} ^{G}$ is
given by:
\begin{multline*}
\sum _{(g)\subset G} \sum _{I\subset  \{1,\dotsc ,n \} } z (g)
\LangRang{c_{I}a_{i}a_{j}
(g)}^{G}\LangRang{(g^{-1})a_{k}a_{l}c_{I^{c}}}^{G}.
\end{multline*}
\end{thm}
From this theorem, one easily derives the PDE version of the WDVV
equations:
\begin{corollary}\label{cor: WDVV PDE formulation}
Let $F_{G}$ be the generating function for the
$G$-Hurwitz-Hodge integrals. Let
\[
F_{ijk} = \frac{\partial ^{3}F_{G}}{\partial x_{i}\partial x_{j}\partial x_{k}}
\]
and let 
\[
g_{ij} = \frac{1}{z_{i}}\delta _{i\jbar}, \quad  g^{ij} = z_{i}\delta _{i\jbar }
\]
where $z_{i}$ is the order of the centralizer of the $i$th conjugacy
class and if $(g) $ is the $j$th conjugacy class then $(g^{-1})$ is
the $\jbar $th conjugacy class. Then the following expression is
symmetric in $\{i,j,n,m \}$:
\[
g_{ij}g_{nm}|G| +\sum _{k,l} F_{ijk}g^{kl}F_{lnm}.
\]
\end{corollary}
The constant term in the above expression corresponds to terms
containing an insertion of the trivial conjugacy class. These terms
occur separately from the derivative terms since our variables only
correspond to non-trivial conjugacy classes.

The proof Theorem~\ref{thm: WDVV sum form} is substantially no
different than the proof of the WDVV equations in orbifold
Gromov-Witten theory given in \cite[\S6.2]{Ab-Gr-Vi-2}. The only
difference is that we are integrating the total Chern class of the
dual Hodge bundle and so we need to check that the Hodge bundle
behaves well on the boundary. Indeed, the Hodge bundle restricted to
the boundary component where two domain curves are glued along a
marked point is equal in K-theory (up to a trivial factor) to the sum
of the Hodge bundles of each factor.

\begin{remark}\label{rem: FG automatically satisfies WDVV}
\emph{We may assume that the series $F_{G}$ given in Theorem~\ref{thm:
main theorem root system formulation} satisfies the WDVV equations,
before we actually prove that it is the generating function for the
$G$-Hurwitz-Hodge integrals.}  The formula in Theorem~\ref{thm: main
theorem root system formulation} for $F_{G}$ was obtained from
$F_{Y}$, the Gromov-Witten potential for the Calabi-Yau threefold
$Y=\GHilb (\C ^{3})$ by a linear change of variables and a
specialization of the quantum parameters (see
\cite{Bryan-Gholampour3}). Since the change of variables transforms
the Poincar\'e pairing on $Y$ to the pairing $g_{ij}$ defined above,
it transforms the WDVV equations for $F_{Y}$ into the WDVV equations
for $F_{G}$. Thus the predicted formula for $F_{G}$ automatically
satisfies the WDVV equations. This is a feature common to all
predictions obtained via the Crepant Resolution Conjecture.
\end{remark}

\section{Computing
$F_{\ZtwoZtwo }$}\label{sec: FZ2Z2} In this section, we fix
$G=\ZtwoZtwo $. For the integral $\LangRang{ c_{1}\dotsb c_{n}
}$ to be non-zero, we must have the monodromy condition
satisfied: the product of the insertions must be trivial. This is
equivalent to
\begin{equation}\label{eqn: monodromy condition for Z2xZ2}
\LangRang{ \zeta ^{n_{1}}_{1} \zeta ^{n_{2}}_{2} \zeta ^{n_{3}}_{3} }
=0\quad \text{ unless }\quad n_{1}\equiv n_{2}\equiv n_{3} \mod 2.
\end{equation}
Consequently, the only non-trivial integrals of length three are
\[
\LangRang{ \zeta _{1}\zeta _{2}\zeta _{3} }=\LangRang{
1\, \zeta _{1}\zeta _{1} }=\LangRang{ 1\, \zeta
_{2}\zeta _{2} }=\LangRang{ 1\, \zeta _{3}\zeta _{3}
} = \frac{1}{4}.
\]
We also note the integrals $\LangRang{ \zeta _{1}^{n_{1}} \zeta
_{2}^{n_{2}} \zeta _{3}^{n_{3}} }$ are symmetric under
permutations of $(\zeta _{1},\zeta _{2},\zeta _{3})$.

\begin{lemma}\label{lem: Z2Z2 determined by WDVV and initial cond}
The integrals $\LangRang{ \zeta _{1}^{n_{1}} \zeta _{2}^{n_{2}}
\zeta _{3}^{n_{3}} }$ are uniquely determined by the
length three integrals, the integrals $\LangRang{ \zeta _{1}^{n}
}$, and the WDVV equations.
\end{lemma}
\textsc{Proof:} We proceed by induction on length (total number of
insertions). The length three integrals start the induction and so we
fix $n>4$ and we assume that all integrals of length less than $n$ are
known. We introduce the notation 
\[
\LL (<n)
\]
to stand for any combination of
integrals of length less than $n$.

Fix $(k_{1},k_{2},k_{3})$ with
\[
k_{1}+k_{2}+k_{3} = n-3
\]
and consider the WDVV relation
\[
\LangRang{ \zeta _{1}^{k_{1}} \zeta _{2}^{k_{2}} \zeta _{3}^{k_{3}}
\,(\zeta _{1}\zeta _{1}|\zeta _{2}\zeta _{2}) }=
\LangRang{ \zeta _{1}^{k_{1}} \zeta _{2}^{k_{2}} \zeta _{3}^{k_{3}}
\,(\zeta _{1}\zeta _{2}|\zeta _{1}\zeta _{2}) }.
\]
Expanding out each side into a sum of products of integrals and
applying the monodromy condition \eqref{eqn: monodromy condition for
Z2xZ2}, we find that there is only one non-zero term of length
$n$. This yields:
\[
\LangRang{ \zeta _{1}^{k_{1}+1} \zeta _{2}^{k_{2}+1}
\zeta _{3}^{k_{3}+1} } = \LL (<n).
\]
Since the integrals $\LangRang{\zeta _{i}^{n}}$ are known by hypothesis, the only unknowns of length $n$ are of the form 
\[
\LangRang{\zeta _{i}^{a}\zeta _{j}^{b}} \text{ where $i\neq j$ and $a+b=n$.}
\]

Now consider the following WDVV relation with $k_{1}+k_{2}=n-3$
\[
\LangRang{\zeta _{1}^{k_{1}}\zeta _{2}^{k_{2}} \,(\zeta _{1}\zeta
_{1}|\zeta _{2}\zeta _{3})} = \LangRang{\zeta _{1}^{k_{1}}\zeta _{2}^{k_{2}} \,(\zeta _{1}\zeta
_{2}|\zeta _{1}\zeta _{3})}.
\]
Expanding out we obtain
\[
\LangRang{\zeta _{1}^{k_{1}+3}\zeta _{2}^{k_{2}}} = 
\LangRang{\zeta _{1}^{k_{1}+1}\zeta _{2}^{k_{2}+2}} +
\LangRang{\zeta _{1}^{k_{1}+1}\zeta _{2}^{k_{2}}\zeta _{3}^{2}} +\LL (<n).  
\]
Solving for $\LangRang{\zeta _{1}^{k_{1}+1}\zeta _{2}^{k_{2}+2}}$,
using the previous equation to write $\LangRang{\zeta
_{1}^{k_{1}+1}\zeta _{2}^{k_{2}}\zeta _{3}^{2}}$ in terms of $\LL
(<n)$, and setting $k_{1}=a-1$ and $k_{2}=b-2$ we get
\[
\LangRang{\zeta _{1}^{a}\zeta _{2}^{b}} = \LangRang{\zeta _{1}^{a+2}\zeta _{2}^{b-2}} +\LL (<n)
\]
for any $a\geq 1$ and $b\geq 2$ with $a+b=n$. By the monodromy
condition \eqref{eqn: monodromy condition for Z2xZ2}, we have that $a$
and $b$ must both be even, so we can use the above equation to
inductively solve for all $\LangRang{\zeta _{1}^{a}\zeta _{2}^{b}}$ in
terms of $\LangRang{\zeta _{1}^{n}}$ and integrals of length less than
$n$ and the lemma is proved.\qed

To prove Proposition~\ref{prop: explicit formula for FZ2Z2}, we now
show that the series in Proposition~\ref{prop: explicit formula for
FZ2Z2} is the unique solution to the WDVV relations which has the
correct cubic terms and the correct specialization $F_{A_{4}}
(x,0,0)$.

Up to symmetry, there are two distinct WDVV relations for the
$\ZtwoZtwo $-Hurwitz-Hodge integrals. In generating function form
(Corollary~\ref{cor: WDVV PDE formulation}), the relations are
\[
F_{121}^{2}+F_{122}^{2}+F_{123}^{2} = F_{111}F_{122}+F_{112}F_{222}+F_{113}F_{322}+\frac{1}{16}
\]
and 
\[
F_{121}F_{133}+F_{122}F_{233}+F_{123}F_{333} = F_{131}F_{123}+F_{132}F_{223}+F_{133}F_{323}.
\]
It is a straight forward but tedious exercise in trigonometry to prove
that the series given in Proposition~\ref{prop: explicit formula for
FZ2Z2} satisfies the above WDVV equations (see Remark~\ref{rem: FG
automatically satisfies WDVV} for a conceptual proof).

To finish the proof of Proposition~\ref{prop: explicit formula for
FZ2Z2}, it remains to check that the formula for $F$ given in
Proposition~\ref{prop: explicit formula for FZ2Z2} has the correct
specializations, namely that
\[
F_{ijk} (0,0,0) = \LangRang{\zeta _{i}\zeta _{j}\zeta _{k}}
\]
and 
\[
F (u,0,0)=\sum _{n\geq 3}\LangRang{\zeta _{1}^{n}}\frac{u^{n}}{n!}.
\]
The first is easy to check and the second is equivalent to the following:
\begin{lemma}
Let $F (x_{1},x_{2},x_{3})$ be the series given in
Proposition~\ref{prop: explicit formula for FZ2Z2}. Then
\[
F_{111} (u,0,0) = \sum _{n=3}^{\infty }\LangRang{\zeta
_{1}^{n}}\frac{u^{n-3}}{(n-3)!}.
\]
\end{lemma}
\textsc{Proof:}
We compute $F_{111} (u,0,0)$ to get
\[
\frac{1}{8}\left(\tan \left(-\frac{u}{4}+\frac{\pi }{4} \right)-\tan
\left(-\frac{u}{4}-\frac{\pi }{4} \right) \right)+\frac{1}{4}\tan
\left(-\frac{u}{2} \right) = \frac{1}{2}\tan (-\frac{u}{2}).
\]
Since for $S= (\zeta _{1}^{n})$ all the monodromies are equal, the
universal cover $\overline{\mathcal{C}}\to \M _{S} (BG)$ is
disconnected and is the union of two copies of the universal double
cover over the hyperelliptic locus $\H _{g}$ where $2g+2=n$.
Consequently, the dual Hodge bundle of $\overline{\mathcal{C}}$ is
\emph{two} copies of the dual Hodge bundle on the hyperelliptic
locus. Denoting both by $\Evee $ we can then write:
\begin{align*}
\sum _{n=3}^{\infty }\LangRang{\zeta _{1}^{n}}\frac{u^{n-3}}{(n-3)!}
&= \sum _{n=3}^{\infty } \frac{u^{n-3}}{(n-3)!} \int _{[\M _{S}
(B\ZtwoZtwo )]}c (\Evee )\\
&= \sum _{g=1}^{\infty }\frac{u^{2g-1}}{(2g-1)!}\int _{[\H _{g}]} c
(\Evee )c (\Evee )\\
&=\sum _{g=1}^{\infty }\frac{u^{2g-1}}{(2g-1)!} \int _{[\H _{g}]}-\lambda _{g}\lambda _{g-1}\\
&=\frac{1}{2}\tan \left(-\frac{u}{2} \right).
\end{align*}
The last equality follows from the computation of Faber-Pandharipande
\cite{Faber-Pandharipande-logarithmic} and the lemma is proved. \qed

This completes the proof of Proposition~\ref{prop: explicit formula
for FZ2Z2}.

\section{Computing $F_{A_{4}}$}\label{sec: FA4}

In this section, we compute the $A_{4}$-Hurwitz-Hodge integrals 
\[
\LangRang{\sigma _{1}^{a_{1}}\sigma _{2}^{a_{2}}\zeta ^{b}}^{A_{4}}
\]
using the WDVV equations along with certain $S_{4}$-Hurwitz-Hodge
integrals that will be computed in \S\ref{sec: computing
S4}. The main technical result of this section is the following.
\begin{prop}\label{prop: FA4 determined by WDVV and specializations}
The generating function $F_{A_{4}} (x_{1},x_{2},x_{3})$ for the
$A_{4}$-Hurwitz-Hodge integrals is uniquely determined by the WDVV
equations, the cubic coefficients of $F_{A_{4}}$, and the
specializations $F_{A_{4}} (x,x,0)$ and $F_{A_{4}} (0,0,x)$.
\end{prop}

We will prove this proposition in \S\ref{subsec: WDVV for
A4}. The specializations that appear in the above proposition can be
expressed in terms of $S_{4}$-Hurwitz-Hodge integrals and $\ZtwoZtwo
$-Hurwitz-Hodge integrals by the following lemma.

\begin{lemma}\label{lem: FA4(x,x,0)=FS4, FA4(0,0,x)=FZ2Z2}
The following equalities hold:
\begin{align*}
3 F_{A_{4}} (0,0,x)&=F_{\ZtwoZtwo } (x,x,x) , \\
2 F_{S_{4}} (0,x,0,0)&=F_{A_{4}} (x,x,0).
\end{align*}
\end{lemma}
\textsc{Proof:} The integrals that appear as coefficients of the
specialization $F_{A_{4}} (0,0,x)$ are $\LangRang{\zeta
^{n}}^{A_{4}}$. These correspond to $A_{4}$ covers whose monodromy
around every branched point is in $\zeta $, the conjugacy class of
disjoint pairs of two cycles. The structure group of such a cover
reduces to $\ZtwoZtwo \subset A_{4}$ and so the integral is given as a
sum of $\ZtwoZtwo $ integrals as follows:
\[
3\LangRang{\zeta ^{n}} ^{A_{4}} = \sum _{n_{1}+n_{2}+n_{3} = n}
\binom{n}{n_{1},n_{2},n_{3}} \LangRang{\zeta
_{1}^{n_{1}}\zeta _{2}^{n_{2}}\zeta _{3}^{n_{3}}}^{\ZtwoZtwo }.
\]
The multinomial coefficient takes into account all the possible
choices of the distribution of the monodromy among the $n$ marked
points. The factor of 3 occurs because the degree of the
map
\[
\M _{0,n} (B\ZtwoZtwo ) \to \M _{0,n} (BA_{4})
\]
is 3. By a similar argument, we derive
\[
2\LangRang{\sigma  ^{n}} ^{S_{4}} = \sum _{n_{1}+n_{2} = n}
\binom{n}{n_{1},n_{2}} \LangRang{\sigma 
_{1}^{n_{1}}\sigma  _{2}^{n_{2}}}^{A_{4}}.
\]
The lemma follows easily.\qed 

\begin{corollary}
The validity of Proposition~\ref{prop: explicit formula for FA4},
which gives our formula for $F_{A_{4}}$, follows from
Proposition~\ref{prop: FA4 determined by WDVV and specializations},
Proposition~\ref{prop: explicit formula for FZ2Z2}, and
Proposition~\ref{prop: FS4(0,x,0,0)}.
\end{corollary}
\textsc{Proof:} By Proposition~\ref{prop: FA4 determined by WDVV and
specializations}, we only need to show that the explicit formula for
$F_{A_{4}} (x_{1},x_{2},x_{3})$ given in Proposition~\ref{prop:
explicit formula for FA4} 
\begin{enumerate}
\item satisfies the WDVV equations, 
\item has the correct cubic terms, and
\item has the correct specializations $F_{A_{4}} (x,x,0)$ and
$F_{A_{4}} (0,0,x)$.
\end{enumerate}
The fact that the predicted formula satisfies the WDVV equations is
once again a tedious but straightforward exercise in trigonometry, or
for a more conceptual proof see Remark~\ref{rem: FG automatically
satisfies WDVV}.

The cubic terms correspond to the three point $A_{4}$-Hurwitz-Hodge
integrals. These are simply counts of $A_{4} $ covers which can be
evaluated using group theory and TQFT methods \cite[section
4]{Dijkgraaf-mirror95}. The non-zero values are given by,
\begin{equation}\label{eqn: three point A4 integrals}
\LangRang{\sigma _{1}\sigma _{2}\zeta }  =1, \quad
\LangRang{\sigma _{1}^{3}}  = \LangRang{\sigma
_{2}^{3}} =\frac{4}{3}, \quad \LangRang{\zeta
^{3}} =\frac{1}{2}.
\end{equation}
It is easy to check that cubic terms of the predicted formula for
$F_{A_{4}}$ agrees with the above values.

Finally, in light of Lemma~\ref{lem: FA4(x,x,0)=FS4,
FA4(0,0,x)=FZ2Z2}, we must check that when we specialize the predicted
formula for $F_{A_{4}} (x_{1},x_{2},x_{3})$ to $F_{A_{4}} (0,0,x)$ and
$F_{A_{4}} (x,x,0)$ we get $\frac{1}{3}F_{\ZtwoZtwo } (x,x,x)$ and
$\frac{1}{2}F_{S_{4}} (0,x,0,0)$ which are determined by
Proposition~\ref{prop: explicit formula for FZ2Z2} and
Proposition~\ref{prop: FS4(0,x,0,0)} respectively. With the use of the
following trigonometric identities:
\begin{align}\label{eqn: trig identities}
\frac{1}{9}\h (3u) &= \h (u) + \h \left(u+\frac{2\pi }{3} \right)+ \h \left(u-\frac{2\pi }{3} \right),\\ \nonumber
\frac{1}{4}\h (2u) &= \h \left(u+\frac{\pi }{2} \right)+ \h \left(u-\frac{\pi }{2} \right),
\end{align}
this is a straightforward check. \qed 

\subsection{The WDVV relations for $A_{4}$-Hurwitz-Hodge
integrals}\label{subsec: WDVV for A4} In this subsection we give the
proof of Proposition~\ref{prop: FA4 determined by WDVV and
specializations}. As before, we use the notation
\[
\LL (<n)
\]
to denote any combination of integrals of length less than $n$. 

We will use induction on the length to prove that the integrals are
determined by the WDVV relations from the length three integrals
(which start the induction) and the integrals (or combinations of
integrals) which occur as the coefficients of the specializations
$F_{A_{4}} (0,0,x)$ and $F_{A_{4}} (x,x,0)$.

The WDVV relations we need are given in the following:

\begin{lem} \label{lem:relations}
Let $n=a_1+a_2+b+3$. We have the following relations among the
$A_{4}$-Hurwitz-Hodge integrals:
\begin{enumerate}
\item [i)] $4\langle \s_1^{a_1+2} \, \s_2^{a_2} \, \z^{b+1}
\rangle  = 4\langle \s_1^{a_1} \, \s_2^{a_2+1} \, \z^{b+2}
\rangle  + \LL(<n),$ \item [ii)]$ 4\langle \s_1^{a_1}
\,\s_2^{a_2+2} \, \z^{b+1} \rangle  = 4\langle \s_1^{a_1+1}
\, \s_2^{a_2} \, \z^{b+2} \rangle  + \LL(<n),$ \item [iii)]$
4\langle \s_1^{a_1+1} \,\s_2^{a_2+1} \, \z^{b+1} \rangle 
 = 4\langle \s_1^{a_1} \, \s_2^{a_2} \, \z^{b+3} \rangle  +
\LL(<n),$ \item [iv)]$ 4\langle \s_1^{a_1+3} \, \s_2^{a_2} \,
\z^{b} \rangle  + 4\langle \s_1^{a_1} \,\s_2^{a_2+3} \,
\z^{b} \rangle  = 8 \langle \s_1^{a_1+1} \, \s_2^{a_2+1} \,
\z^{b+1} \rangle  +\LL(<n).$
\end{enumerate}
\end{lem}
\textsc{Proof:} We prove the above relations using the following WDVV
relations which are expressed using the notation of Theorem~\ref{thm:
WDVV sum form}.
\begin{enumerate}
\item[i)]$\langle \s_1^{a_1} \,\s_2^{a_2} \, \z^{b}\, (\s_1\,\z\, |\, \s_1\, \z) \rangle =
\langle \s_1^{a_1} \,\s_2^{a_2} \, \z^{b}\, (\s_1\,\s_1\, |\, \z \, \z) \rangle$,
\item[ii)]$\langle \s_1^{a_1} \,\s_2^{a_2} \, \z^{b}\, (\s_2\,\z\, |\, \s_2\, \z) \rangle =
\langle \s_1^{a_1} \,\s_2^{a_2} \, \z^{b}\, (\s_2\,\s_2\, |\, \z\, \z) \rangle$,
\item[iii)]
$\langle \s_1^{a_1} \,\s_2^{a_2} \, \z^{b}\, (\s_1\,\z\, |\, \s_2\, \z) \rangle = \langle \s_1^{a_1}
\,\s_2^{a_2} \, \z^{b}\, (\s_1\,\s_2\, |\, \z\, \z) \rangle$,
\item[iv)]$\langle \s_1^{a_1} \,\s_2^{a_2} \, \z^{b}\, (\s_1\,\s_1\, |\, \s_2\, \s_2) \rangle =
 \langle \s_1^{a_1} \,\s_2^{a_2} \, \z^{b}\, (\s_1\,\s_2\, |\, \s_1\,
\s_2) \rangle$.
\end{enumerate}
After expanding i) for $a_1+a_2+b>0$, the resulting equation is
\begin{align*}
\sum
\binom{a_1}{a'_1}\binom{a_2}{a'_2}\binom{b}{b'}\leftbig{\{}{20}
&4 \langle \s_1^{a'_1+1} \,\s_2^{a'_2} \,
\z^{b'+2} \rangle \langle \s_1^{a''_1+1} \,\s_2^{a''_2} \, \z^{b''+2} \rangle+\\
&3\langle \s_1^{a'_1+2} \,\s_2^{a'_2} \, \z^{b'+1} \rangle \langle
\s_1^{a''_1+1} \,\s_2^{a''_2+1} \, \z^{b''+1} \rangle+\\&3\langle
\s_1^{a'_1+1} \,\s_2^{a'_2+1} \, \z^{b'+1} \rangle \langle
\s_1^{a''_1+2} \,\s_2^{a''_2} \, \z^{b''+1} \rangle\rightbig{\}}{20}=
\end{align*}
\begin{align*}
\sum \binom{a_1}{a'_1}\binom{a_2}{a'_2}\binom{b}{b'}\leftbig{\{}{20}&4 \langle \s_1^{a'_1+2} \,\s_2^{a'_2} \,
\z^{b'+1} \rangle \langle \s_1^{a''_1} \,\s_2^{a''_2} \, \z^{b''+3} \rangle+\\
&3\langle \s_1^{a'_1+3} \,\s_2^{a'_2} \, \z^{b'} \rangle \langle \s_1^{a''_1} \,\s_2^{a''_2+1} \, \z^{b''+2}
\rangle+\\&3\langle \s_1^{a'_1+2} \,\s_2^{a'_2+1} \, \z^{b'} \rangle \langle \s_1^{a''_1+1} \,\s_2^{a''_2} \,
\z^{b''+2} \rangle \rightbig{\}}{20}
\end{align*}
where the sums are over $a'_1+a''_1=a_1$, $a'_2+a''_2=a_2$, and
$b'+b''=b$.

All the integrals of length $n$ in the above expression are multiplied
by a integrals of length 3 which are given in equation~\eqref{eqn:
three point A4 integrals}. All other terms contain only terms of
length less than $n$. Substituting the values of the length three
integrals, we obtain the relation i). The proof of relations ii),
iii), and iv) is similar. \qed

\begin{rem} \label{rem:monodoromy}
One can see easily by monodromy considerations that $$\langle
\s_1^{a_1} \, \s_2^{a_2} \, \z^{b} \rangle $$ is nonzero only when
$a_1 \equiv a_2 \; (\text{mod 3}).$ Note also that the above integral
is symmetric in $a_1$ and $a_2$ due to the fact that $A_4$ has a
nontrivial outer automorphism which exchanges $\sigma _{1}$ and
$\sigma _{2}$.
\end{rem}

We will now use the relations i)--iv) in Lemma~\ref{lem:relations} to
show that all the $A_{4} $ integrals can be inductively recovered from
the length three integrals, the integrals $\LangRang{\zeta
^{n}}^{A_{4}}$ (which are the coefficients of $F_{A_{4}} (0,0,x)$),
and the integrals $\LangRang{\sigma ^{n}}^{S_{4}}$ (which by virtue
of Lemma~\ref{lem: FA4(x,x,0)=FS4, FA4(0,0,x)=FZ2Z2} are the
coefficients of $F_{A_{4}} (x,x,0)$).

The following relations are direct consequences of relations i) and
iii) in Lemma~\ref{lem:relations}:
\begin{align*}
\langle \s_1^{k} \, \z^{b} \rangle &=\langle \s_1^{k-2} \,
\s_2 \, \z^{b+1} \rangle +\LL(<k+b) \quad \quad \quad \quad \quad \quad b>0, \; k>2, \nonumber \\
\langle \s_1^{a_1} \, \s_2^{a_2} \, \z^{b} \rangle & =\langle
\s_1^{a_1-1} \, \s_2^{a_2-1} \, \z^{b+2}
\rangle +\LL(<a_1+a_2+b) \quad \quad a_1,\, a_2, \, b>0.
\end{align*}

The following lemma is proven readily by a direct repeated
application of two relations above:

\begin{lem} \label{lem:case b>0}
Assume that $b>0$. We have the following relation:
$$\langle \s_1^{a_1} \, \s_2^{a_2} \, \z^{b} \rangle   =\langle \z^{a_1+a_2+b}
\rangle +\LL(<a_1+a_2+b).$$
\end{lem} \qed

In the next lemma we deal with the case $b=0$:

\begin{lem} \label{lem:case b=0}
Let $n=a_1+a_2$. All the integrals of the form $\langle \s_1^{a_1} \,
\s_2^{a_2} \rangle $ can be written in terms of $\langle \s^{n}
\rangle^{S_4}$, length $n$ integrals having at least one $\zeta $
insertion, and $\LL(<n)$.
\end{lem}
\textsc{Proof:} Let $0 \le k \le 2$ be so that $n \equiv k \;
(\text{mod 3}).$ One can see easily that (see
Remark~\ref{rem:monodoromy}):

\begin{align*}
\{ \langle \s_1^{a_1} \, \s_2^{a_2} & \rangle  |\; a_1+a_2=n \} \\
&= \{ \langle \s_1^{n-k} \, \s_2^{k} \rangle ,\; \langle
\s_1^{n-k-3} \, \s_2^{k+3} \rangle , \;\dots , \; \langle
\s_1^{k} \, \s_2^{n-k} \rangle \}.
\end{align*}
This is a set of $l+1$ elements, where 
\[
l=(n-2k)/3.
\]
For simplicity we write the elements of this set in the same order
form left to right by $x_0,\, x_2,\dots,\, x_{l}$. Note that by
symmetry of the integrals (see Remark~\ref{rem:monodoromy}) we have
\begin{equation} \label{equ:symmetry} 
x_i=x_{l-i}.
\end{equation}

Applying relation iv) in Lemma~\ref{lem:relations} $l$ times, we
get the following set of relations: \begin{align}
\label{equ:relations, b=0}
x_0+x_1&=2\langle \s_1^{n-k-2} \, \s_2^{k+1} \, \z \rangle +\LL(<n) \nonumber\\
x_1+x_2&=2\langle \s_1^{n-k-5} \, \s_2^{k+4} \, \z \rangle +\LL(<n) \nonumber\\
& \; \; \vdots \nonumber \\
x_{l-1}+x_l&=2\langle \s_1^{k+1} \, \s_2^{n-k-2} \, \z
\rangle +\LL(<n).
\end{align}
Now we consider two cases: 
\begin{enumerate} 
\item [i)] $n$ is odd: \\
We see that $l$ (the number of relations in (\ref{equ:relations,
b=0})) is odd as well, and then (\ref{equ:symmetry}) and
(\ref{equ:relations, b=0}) boil down into the following set of
independent equations among $x_i$'s:
\begin{align*}
x_0+x_1&=2\langle \s_1^{n-k-2} \, \s_2^{k+1} \, \z \rangle +\LL(<n) \\
x_1+x_2&=2\langle \s_1^{n-k-5} \, \s_2^{k+4} \, \z \rangle +\LL(<n) \\
& \; \; \vdots  \\
x_{p}+x_p&=2\langle \s_1^{n-k-2-3p} \, \s_2^{k+1+3p} \, \z
\rangle +\LL(<n),
\end{align*}
where $p=(l+1)/2$. From these equations we get
\begin{align*}
x_{l-i}=x_i &=-(-1)^{p-i}\langle \s_1^{n-k-2-3p} \, \s_2^{k+1+3p}
\, \z \rangle \\& \quad +2\sum_{j=i}^p (-1)^{j-i}\langle
\s_1^{n-k-2-3j} \, \s_2^{k+1+3j} \, \z \rangle +\LL(<n),
\end{align*}
and the Lemma is proven in this case.

\item [ii)] $n$ is even: \\
In this case $l$ is even. One can see that one of the relations in
(\ref{equ:symmetry}) is redundant. However, since the coefficients of
$F_{A_{4}} (x,x,0) = 2F_{S_{4}} (0,x,0,0)$ are known by hypothesis, we
get the following extra relation:
\begin{equation} \label{equ:<3^n>}
\binom{n}{n-k}x_0 \, + \,  \binom{n}{n-k-3}x_1 \, + \dots  +\,
\binom{n}{k}x_l \;=\;2\langle \s^{n} \rangle^{S_4}.
\end{equation}
Adding (\ref{equ:<3^n>}) to (\ref{equ:relations, b=0}), we get a
system of $l+1$ equations and $l+1$ unknowns (we call $x_i$'s
unknowns). One can see that the determinant of the matrix of
coefficients is given by (note that n is even)
\[
\binom{n}{n-k}\,-\,  \binom{n}{n-k-3}\,+\dots +\, \binom{n}{k},
\]
which is nonzero by \cite[Lemma~A.6]{Bryan-Graber-Pandharipande}. This
means that we can express $x_i$'s in terms of the right hand side of
the system of $l+1$ equations, and the Lemma is proven in this
case.\end{enumerate} \qed

We may now prove Proposition~\ref{prop: FA4 determined by WDVV and
specializations} in the following equivalent form:

\begin{prop}
The $A_{4}$-Hurwitz-Hodge integrals are uniquely determined by the
WDVV equations, the length three integrals, and the integrals $\langle
\z^{m} \rangle ^{A_{4}} $ and $\langle \s^{m} \rangle^{S_4}$.
\end{prop}
\textsc{Proof:} We use induction on the length $n$ of the
integrals. The length three integrals are known by hypothesis.

Let $\langle \s_1^{a_1} \, \s_2^{a_2} \, \z^{b} \rangle ^{A_{4}}$ be
an arbitrary integral of length $n>3$. By Lemmas \ref{lem:case b>0}
and \ref{lem:case b=0}, we can write this integral in terms of
$\LangRang{\zeta ^{n}}^{A_{4}}$, $\LangRang{\sigma ^{n}}^{S_4}$, and
$\LL(<n)$. Both $\langle \z^{n} \rangle^{A_{4}} $ and $\langle \s^{n}
\rangle^{S_4}$ are known by the assumption, and $\LL(<n)$ is also
known by the induction hypothesis. Therefore $\langle \s_1^{a_1} \,
\s_2^{a_2} \, \z^{b} \rangle^{A_{4}} $ is determined, and the
Proposition is proven. \qed

\section{Computing $S_{4}$-Hurwitz-Hodge integrals}\label{sec: computing S4}

In this section we prove the following proposition, which is needed to
complete the proof of Proposition~\ref{prop: explicit formula for
FA4}.

\begin{prop}\label{prop: FS4(0,x,0,0)}
Let $F_{S_{4}} (x_{1},x_{2},x_{3},x_{4})$ be the generating function
for the $S_{4}$-Hurwitz-Hodge integrals. Then
\[
F_{S_{4}} (0,x,0,0) = \frac{1}{8}\h \left(\frac{2x}{\sqrt{3}} -
\frac{2\pi }{3} \right) +2 \h \left(\frac{x}{\sqrt{3}} - \frac{\pi
}{3} \right).
\]
\end{prop}
This follows immediately from the identity in equations~\eqref{eqn:
trig identities} and the following:
\begin{theorem}\label{thm: FS4(0,u,0,v)} 
\[
F_{S_{4}} (0,u,0,v) = \frac{1}{2}K (2u,v) +K (-u,v)
\]
where
\begin{align*}
K (u,v) =& \quad \h \left(\frac{u}{\sqrt{3}}+\frac{v}{2} -\frac{5\pi }{6} \right) +2 \h \left(\frac{u}{\sqrt{3}}-\frac{\pi }{3} \right)\\
&+\h \left(\frac{u}{\sqrt{3}}-\frac{v}{2} +\frac{\pi }{6} \right) +2 \h \left(\frac{u}{\sqrt{2}}+\frac{\pi }{2} \right)\\
&+\frac{2}{3}\h \left(\frac{v}{2}+\frac{3\pi }{2} \right).
\end{align*}
\end{theorem}

To prove Theorem~\ref{thm: FS4(0,u,0,v)}, the basic strategy is once
again to use WDVV to inductively determine the integrals. We will not
determine all the $S_{4}$-Hurwitz-Hodge integrals but we will have to
determine a certain set of integrals that have a small number of $\tau
$ and $\rho $ insertions.

We use the following generating functions. For convenience, we define
unstable integrals (those have fewer than three insertions) to be
zero.

\begin{itemize}
\item $\displaystyle T(u,v)= F_{S_{4}} (0,u,0,v)=\sum_{a,b =
0}^{\infty } \langle \s^a \; \z^b \rangle\,
\frac{u^a}{a!}\frac{v^b}{b!},$
\item $\displaystyle X_a(u)= \sum_{n = 0}^{\infty } \langle \s^a \,
\z^{n} \rangle \frac{u^{n}}{n!},$
\item $\displaystyle Y_b(u)= \sum_{n = 0}^{\infty } \langle \s^n \,
\z^{b} \rangle \frac{u^{n+b-3}}{(n+b-3)!},$
\item $\displaystyle B(u)= \sum_{n =0}^{\infty } \langle \T^2\, \s^n
\rangle \frac{u^{n-1}}{(n-1)!},$
\item $\displaystyle C(u)= \sum_{n = 0}^{\infty } \langle \T\; \s^n \,
\R \; \z \rangle\, \frac{u^{n}}{n!},$
\item $\displaystyle D(u)= \sum_{n =0}^{\infty } \langle \T \; \s^n \,
\R \rangle \frac{u^{n-1}}{(n-1)!}.$
\end{itemize}

The length three integrals can be evaluated using group theory and
TQFT methods \cite[section 4]{Dijkgraaf-mirror95}. The non-zero values
are given by:
\begin{align*}
\LangRang{\sigma ^{2}\zeta }=\LangRang{\tau ^{2}\sigma }&=\LangRang{\tau \rho \sigma } =\LangRang{\rho ^{2}\sigma } =1,\\
\LangRang{\zeta ^{3}} = \LangRang{\rho ^{2}\zeta } = \LangRang{\tau ^{2}\zeta } &= \LangRang{\tau ^{2}1} = \LangRang{\rho ^{2}1}=\frac{1}{4},\\
\LangRang{\sigma ^{3}} = \frac{4}{3},\quad \LangRang{\tau \rho \zeta }
=\frac{1}{2},&\quad \LangRang{\sigma ^{2}1}=\frac{1}{3}, \quad
\LangRang{\zeta ^{2}1} = \frac{1}{8}.
\end{align*}

We now determine the series $X_{0}$, $X_{1}$, and $B$.

\begin{lemma}\label{lem: the series X0, X1, and B}
The series $X_{0} (u)$, $X_{1} (u)$, and $B (u)$ are given by
\begin{align*}
X_{0} (u) =& \frac{1}{6} \h \left(\frac{3u}{2} - \frac{\pi }{2} \right) +\frac{1}{2}\h \left(\frac{u}{2} +\frac{\pi }{2} \right) +\frac{1}{4}\h \left(u \right),\\
X_{1} (u) =& 0,\\
B (u) =& \frac{1}{\sqrt{3}}\tan(\frac{-u}{\sqrt{12}}+\frac{\pi}{3}).&
\end{align*}
\end{lemma}
\textsc{Proof:} The proof of the first formula follows immediately
from Proposition~\ref{prop: explicit formula for FZ2Z2} and the
formula
\[
X_{0} (u) = \frac{1}{6}F_{\ZtwoZtwo } (u,u,u).
\] 
The proof of the above formula is almost identical to the proof of the
first formula in Lemma~\ref{lem: FA4(x,x,0)=FS4, FA4(0,0,x)=FZ2Z2};
the only difference being that the 3 (the index of $\ZtwoZtwo $ in
$A_{4}$) is replaced with 6 (the index of $\ZtwoZtwo $ in
$S_{4}$). 

The Lemma's second formula is a consequence of
\[
\LangRang{\sigma \zeta ^{n}} = 0
\]
which follows from monodromy considerations (see
Remark~\ref{rem:monodoromy}).

To prove the Lemma's third formula, we need to show that the
generating function $B(u)$ has the same coefficients, up to an
alternating sign, as the generating function with the same name
computed in \cite[Proposition~A.3]{Bryan-Graber-Pandharipande}. The
coefficients of the generating function in
\cite{Bryan-Graber-Pandharipande} are the $\lambda _{g}$ integrals
over the Hurwitz locus of curves admitting a \emph{degree three} cover
of $\P$ with 2 ordinary ramification points and $g+1$ double
ramifications. The identification of this integral with the integral
$(-1)^{g}\LangRang{\tau ^{2}\sigma ^{g+1}}^{S_{4}}$ is because the
only chance for a degree four cover to contribute to $\LangRang{\tau
^{2}\sigma ^{g+1}}^{S_{4}}$ is if the cover is
\emph{disconnected}. Indeed, if it is not disconnected, then the genus
is $g-1$ and hence the $g$th Chern class of $\Evee $ is zero. The sign
is because we work with the dual Hodge bundle instead of the Hodge
bundle.\qed

To determine $T (u,v)$, it clearly suffices to determine $X_{a} (u)$
for all $a$. We will do this using an induction on $a$. The following
lemma provides the basic relations that are needed in the induction.

\begin{lemma}\label{lem: WDVV relations for S4 integrals}
The series $T (u,v)$, $X_{a} (u)$, $Y_{0} (u)$, $Y_{1} (u)$, $Y_{2} (u)$,
and $Y_{3} (u)$ satisfy the following relations.
\begin{enumerate}
\item [(i)] $ 3T^2_{uuv}+8T^2_{uvv}-3T_{uuu}T_{uvv}-8T_{uuv}T_{vvv}+1=0,$ \vspace{.1in}
\item [(ii)] $3\left(X_{2}' \right)^{2}-8X_{2}'X_{0}''' -1
=0,$\vspace{.1in}
\item [(iii)]$(6Y_2B+4Y_3+1)(3Y_0B+2Y_1-4B^2+2)=2(3Y_1B+2Y_2-B)^2,$\vspace{.1in}
\item [(iv)] $(6X_{2}'-8X_{0}''')X_{a+2}'-3aX_{2}''X_{a+2} = G
(X_{0},X_{1},\dotsc ,X_{a+1}),$
\end{enumerate}
where in item (iv) $a>0$, and $G (X_{0},X_{1},\dots
,X_{a+1})$ is a function of the series $X_{k} (u)$ for $0\leq k\leq
a+1$.
\end{lemma}
\textsc{Proof:} These relations are all consequences of the WDVV
equations. For fixed $a$ and $b$, consider the WDVV relation
\[
\LangRang{\sigma ^{a}\zeta ^{b} (\sigma \zeta |\sigma \zeta )} =
\LangRang{\sigma ^{a}\zeta ^{b} (\sigma \sigma |\zeta \zeta )}.
\]
It is given by 
\begin{multline*}
0= \delta _{0,a}\delta _{0,b}+\sum _{\begin{smallmatrix} a_{1}+a_{2}=a\\
b_{1}+b_{2}=b \end{smallmatrix}} \binom{a}{a_{1}}\binom{b}{b_{1}}\cdot \\
\leftbig{\{ }{20} 3\LangRang{\sigma ^{a_{1}+2}\zeta ^{b_{1}+1}}
\LangRang{\sigma ^{a_{2}+2}\zeta ^{b_{2}+1}} + 8\LangRang{\sigma
^{a_{1}+1}\zeta ^{b_{1}+2}}\LangRang{\alpha
^{a_{2}+2}\zeta ^{b_{2}+2}} \\
- 3 \LangRang{\sigma ^{a_{1}+3}\zeta ^{b_{1}}}\LangRang{\alpha
^{a_{2}+1}\zeta ^{b_{2}+2}} -8\LangRang{\sigma ^{a_{1}+2}\zeta
^{b_{1}+1}}\LangRang{\sigma ^{a_{2}}\zeta ^{b_{2}+3}}
\rightbig{\}}{20}
\end{multline*}
Multiplying the above equation by $\frac{u^{a}}{a!}\frac{u^{b}}{b!}$
and then summing over all of $a$ and $b$ yields the PDE in (i). Fixing
$a=0$, multiplying by $\frac{u^{b}}{b!}$, and summing over $b$ yields
(ii). To prove (iv), we fix $a>0$, we multiply the above equation by
$\frac{u^{b}}{b!}$, we move to the right hand side of the equation all
the terms containing only integrals with fewer than $a+2$ insertions
of $\sigma $, and then we sum over $b$.

The relation (iii) follows easily from the equations:
\begin{align*}
 6Y_2\,B+4Y_3+1&=8C^2,\\
 3Y_1\,B+2Y_2 &=4CD+B,\\
 3Y_0\,B+2Y_1+2&=4B^2+4D^2
\end{align*}
which are derived in a fashion similar to the above from the WDVV
relations:
\begin{align*}
\LangRang{\sigma ^{a} (\tau \tau |\zeta \zeta )} &= \LangRang{\sigma ^{a} (\tau \zeta |\tau \zeta )},\\
\LangRang{\sigma ^{a} (\tau \tau |\sigma  \zeta )} &= \LangRang{\sigma ^{a} (\tau \sigma  |\tau \zeta )},\\
\LangRang{\sigma ^{a} (\tau \tau |\sigma \sigma )} &= \LangRang{\sigma ^{a} (\tau \sigma |\tau \sigma )}.
\end{align*}
In the derivation, we use the following crucial fact:
\begin{equation*}
\LangRang{\tau ^{2}\zeta \sigma ^{a}} = \begin{cases}
0&\text{if $a>0$,}\\
\frac{1}{4} &\text{if $a=0$.}
\end{cases}
\end{equation*}
This follows from the fact that a four fold cover of $\P $ with $3$
branched points of monodromies $\tau $, $\tau $, $\zeta $, and $a$
branch points of monodromy $\sigma $ is connected and of genus
$a-1$. Consequently $c_{a} (\Evee )=0$, namely the $a$th Chern
class of the dual Hodge bundle vanishes, and so the corresponding
Hurwitz-Hodge integral is zero. \qed 

\begin{prop}\label{prop: Xa is uniquely determined by ODE and initial}
The functions $X_{a} (u)$ are uniquely determined for all $a$ by the
series $B (u)$, $X_{0} (u)$, $X_{1} (u)$, the length three integrals,
and the relations (ii), (iii), and (iv) in Lemma~\ref{lem: WDVV
relations for S4 integrals}.
\end{prop}
\textsc{Proof:} The relation (ii) of Lemma~\ref{lem: WDVV relations
for S4 integrals} is a quadratic equation for $X_{2}' (u)$ whose
solution is fixed by the the condition
\[
X_{2}' (0) = \LangRang{\sigma ^{2}\zeta } =1.
\]
Since the constant term of $X_{2} (u)$ is an unstable integral (and
hence zero by convention), $X_{2} (u)$ is then uniquely determined.

We now proceed to determine $X_{a} (u)$ by induction on $a$. We assume
that $X_{k} (u)$ is known for all $k<a+2$ and we need to show that we
can determine $X_{a+2} (u)$. Since $X_{0}$, $X_{1}$, and $X_{2}$ are
known, we may assume that $a>0$. Then relation (iv) in Lemma~\ref{lem:
WDVV relations for S4 integrals} is a first order ODE for $X_{a+2}
(u)$. Since the coefficient of $X'_{a+2}$ in the ODE is an invertible
series, the ODE has a solution which is uniquely determined by
specifying $X_{a+2} (0)$. 

Now $X_{a+2} (0)$ is equal to the coefficient of $u^{a-1}$ in the
series $Y_{0} (u)$, and so we need to determine this
coefficient. Since the series $X_{0},\dotsc ,X_{a+1}$ are known by the
induction hypothesis, we know the coefficients of $Y_{b} (u)$ up to
the $u^{a+b-2}$ term. By examining the $u^{a-1}$ term in the relation
(iii) from Lemma~\ref{lem: WDVV relations for S4 integrals}, we find
that the only unknown is the coefficient of $u^{a-1}$ of $Y_{0}$ which
appears exactly once with a non-zero coefficient. Hence we can
uniquely solve for this coefficient which provides the initial
condition which uniquely determines $X_{a+2} (u)$. \qed 

Thus to complete the proof of Theorem~\ref{thm: FS4(0,u,0,v)}, we must
show that the formula for $T (u,v)$ in the theorem, yields series
$X_{a} (u)$  
\begin{enumerate}
\item which predict the correct length three integrals,
\item which agree with the formulas for $X_{0}$ and $X_{1}$ given in
Lemma~\ref{lem: the series X0, X1, and B},
\item and which are solutions to the relations $(ii)$, $(iii)$, and
$(iv)$ of Lemma~\ref{lem: WDVV relations for S4 integrals}.
\end{enumerate}
The first two are straightforward checks. The compatibility of the
solution with relation (iii) is also a straightforward check. The
compatibility with relations (ii) and (iv) is equivalent to the
formula for $T$ satisfying the PDE (i). This compatibility can be
checked directly (with Maple for example), but there is a more
conceptual proof, along the lines of Remark~\ref{rem: FG automatically
satisfies WDVV} which we outline below.

The formula for $T (u,v)$ can be derived from the formula for
$F_{A_{4}} (x_{1},x_{2},x_{3})$ via the relation
\[
T (u,v) = F_{S_{4}} (0,u,0,v) = \frac{1}{2}F_{A_{4}} (u,u,v)
\]
which is proved by the same method as the proof of Lemma~\ref{lem:
FA4(x,x,0)=FS4, FA4(0,0,x)=FZ2Z2}.  The fact that $T (u,v)$ satisfies
(i) can be seen to be a consequence of the fact that $F_{A_{4}}$
satisfies the $A_{4}$ WDVV equations, a fact which we proved in
\S\ref{sec: FA4}. Recall that our formula for $F_{A_{4}} $ was
derived via the Crepant Resolution Conjecture from the Gromov-Witten
potential of the crepant resolution $Y\to \C ^{3}/A_{4}$ given by the
$A_{4}$ Hilbert scheme. The derived formula for $F_{A_{4}}$ thus
automatically satisfies the WDVV equations since the Gromov-Witten
potential of $Y$ satisfies the WDVV equations and the Crepant
Resolution Conjecture is compatible with the WDVV equations. Thus the
fact that $T$ satisfies relations (i) in Lemma~\ref{lem: WDVV
relations for S4 integrals} is ultimately a consequence of the fact
that the Gromov-Witten potential of $Y$ (which was computed in
\cite{Bryan-Gholampour3}) satisfies the WDVV equations.\qed

This completes the proof of Theorem~\ref{thm: FS4(0,u,0,v)} and
consequently it completes the proof of Proposition~\ref{prop: explicit
formula for FA4} and hence it completes the proof of our main
result, Theorem~\ref{thm: main theorem root system formulation}.\qed

\section{The Relationship with orbifold Gromov-Witten
theory}\label{sec: orbifold GW theory} Let $G$ (which is $A_{4}$ or
$\ZtwoZtwo $) act on $\C^3$ by the representation obtained from the
embedding $G\subset SO (3)\subset SU (3)$. Let $\X $ be the orbifold
given by the quotient:
\[
\X =[\C^3/G].
\]

The orbifold cohomology $H^{*}_{\orb} (\X )$ has a canonical basis
labelled by conjugacy classes of $G$. Consequently, the insertions for
the orbifold Gromov-Witten invariants of $\X $ are conjugacy classes,
and the genus zero Gromov-Witten invariants of $\X $ take the form
$\LangRang{ c_{1}\dotsb c_{n} }^{\X }$.

In this section we prove:
\begin{prop}\label{prop: orbifold GW invs = Hurwitz-Hodge ints}
The genus zero orbifold Gromov-Witten invariants\footnote{Strictly
speaking, some of the orbifold Gromov-Witten invariants of $\X $ are
not well defined because the corresponding moduli space of twisted
stable maps is non-compact. In these cases, we \emph{define} the
invariants by localization with respect to the $\C^{*}$ action on $\X
$.} of $\X $ are given by the $G$-Hurwitz-Hodge integrals, namely
\[
\LangRang{ c_{1}\dotsb c_{n} }^{\X }=\LangRang{
c_{1}\dotsb c_{n} }^{G}
\]
for any $n$-tuple of non-trivial conjugacy classes.  Consequently, the
(non-classical part of the) genus zero Gromov-Witten potential of $\X
$ is equal to the generating function of the $G$-Hurwitz-Hodge
integrals:
\[
F_{\X } (x_{1},x_{2},x_{3}) = F_{G} (x_{1},x_{2},x_{3}).
\]
\end{prop}
\textsc{Proof:} By definition, the orbifold invariants are given by
\[
\LangRang{ c_{1}\dotsb c_{k} }^{\X } = \int _{[\M _{0,k}
(\X )]^{vir}} \ev_{1}^{*} (c_{1})\cup \dotsb \cup \ev_{k}^{*} (c_{k}).
\]
Using virtual localization with respect to the $\C ^{*}$ action on $\X
$ , we can express the above integral in terms of an integral over the
$\C ^{*}$ fixed locus of $\M _{0,k} (\X )$ which is $\M _{0,k} (BG)$:
\begin{align*}
\LangRang{ c_{1}\dotsb c_{k} }^{\X }&=\int _{[\M _{0,k}
(BG)]} \ev_{1}^{*} (c_{1})\cup \dotsb \cup \ev_{k}^{*} (c_{k})\cup e
(N^{vir})^{-1}\\
&=\int _{[\M _{S} (BG)]}e (-N^{vir})
\end{align*}
where $S= (c_{1}\dotsb c_{k})$ and $N^{vir}$ is the virtual normal
bundle of $\M _{S} (BG)$ in $\M _{0,k} (\X )$ regarded as an element
in K-theory. So to prove the Proposition, we need to show that 
\begin{equation}\label{eqn: int of Nvir = int of hodge}
\int _{[\M _{S } (BG)]} e (-N^{vir}) = \int _{[\M _{S } (BG)]}c
(R^{1}\overline{\pi }_{*}\mathcal{O}_{\overline{C}}).
\end{equation}
    
Let 
\[
V\to BG
\]
be the bundle given by the three dimensional representation of $G$
induced from the embedding $G\subset SO (3)$. By the standard argument in Gromov-Witten theory, the virtual normal bundle is given by
\[
N^{vir} = -R^{\bullet }\pi _{*}f^{*}V
\]
where $\pi:\mathcal{C}\to \M _{0,k} (BG) $ is the universal curve and
$f:\mathcal{C}\to BG$ is the universal map.

Let $H\subset G$ be the subgroup $\Z _{3}\subset A_{4}$ or $\{0
\}\subset \ZtwoZtwo $ respectively. Then the action of $G$ on the
coset space $G/H$ is the usual permutation action of $G $ on the set
of four elements. Consequently, we can construct the universal degree
4 cover $p:\overline{\mathcal{C}}\to \mathcal{C}$ by pulling back the
map 
\[
i:BH\to BG
\]
via $f$. That is, we have the following diagram:
\[
\begin{diagram}
&\overline{\mathcal{C}} &\rTo{\overline{f}}&BH\\
&\dTo_{p} &&\dTo^{i}\\
\pimapA &\mathcal{C}&\rTo{f}&BG\\
&\dTo_{\pi }&&\\
&\M _{0,k} (BG)&&
\end{diagram}
\]


We now compute in K-theory:
\begin{align*}
R^{\bullet }\overline{\pi }_{*}\mathcal{O}_{\overline{\mathcal{C}}}&=
R^{\bullet }\pi _{*}\left(p_{*}\overline{f}^{*}\mathcal{O}_{BH}
\right)\\
&=R^{\bullet }\pi _{*} (f^{*}i_{*}\mathcal{O}_{BH})\\
&=R^{\bullet }\pi _{*}f^{*} (V\oplus \mathcal{O})\\
&=-R^{1}\pi _{*}f^{*}V + \pi _{*}\mathcal{O}_{\mathcal{C}}.
\end{align*}
The equality on the top line uses the fact that $p$ is finite. The
equality on the second line uses the fact that the commutative square
in the diagram is Cartesian. The equality on the third line uses the
fact that the $G$ representation induced by the trivial representation
is $V$ plus the trivial representation. The equality on the fourth
line uses the fact that $\pi :\mathcal{C}\to \M _{0,k} (BG)$ is genus
0.

Finally, we apply the total Chern class to both sides of the above
equality and integrate over $\M _{S} (BG)$. The $\C ^{*}$ equivariant
Euler class is the same as the total Chern class with the appropriate
power of the equivariant parameter appearing in each degree. Since the
virtual complex dimension of $\M _{0,k} (\X )$ is $k$, and the
insertions $c_{i}$ are all from $H^{2}_{\orb} (\X )$, the integral is
degree 0 in $t$ the equivariant parameter and is hence independent of
$t$. Equation~\eqref{eqn: int of Nvir = int of hodge} is proved. \qed

\bibliography{mainbiblio}
\bibliographystyle{plain}

\end{document}